\documentclass[12pt]{article}
\usepackage{amsmath}
\usepackage{amssymb}
\usepackage{graphicx}
\usepackage[cp1251]{inputenc}
\usepackage[russian]{babel}
\newcommand{\il}[2]{\int\limits_{#1}^{#2}}

\newcommand{\ph}{\phantom{a}}
\newcommand{\phh}{\phantom{aaa}}

\newcommand{\sist}[2]{\left\{
\begin{array}{l}
{#1}\\
\ph\\
{#2}
\end{array}
\right.}

\begin{document}

\vskip 20pt

MSC 34C10

\vskip 20pt

\centerline{\bf Oscillatory and non oscillatory criteria for linear}
 \centerline{\bf four dimensional hamiltonian systems }

\vskip 20 pt

\centerline{\bf G. A. Grigorian}
\centerline{\it Institute  of Mathematics NAS of Armenia}
\centerline{\it E -mail: mathphys2@instmath.sci.am}
\vskip 20 pt

\noindent
Abstract.  The Riccati equation method is used for study the oscillatory and non oscillatory behavior of solutions of linear four dimensional hamiltonian systems.  An oscillatory and three  non oscillatory criteria are proved. On examples the obtained results are compared with some well known  ones.
\vskip 20 pt

Key words: Riccati equation, oscillation, non  oscillation, conjoined (prepared, preferred) solution, Liuville's formula.

\vskip 20 pt

{\bf 1. Introduction.} Let $A(t) \equiv \bigl(a_{jk}(t)\bigr)_{j,k=1}^2,\ph B(t)\equiv \bigl(b_{jk}(t)\bigr)_{j,k=1}^2, \ph C(t)\equiv \ph \bigl(c_{jk}(t)\bigr)_{j,k=1}^2, \linebreak t\ge t_0$, be complex valued continuous matrix functions on $[t_0;+\infty)$ and let $B(t)$ and $C(t)$ be Hermitian, i.e., $B(t) = B^*(t), \ph C(t) = C^*(t), \ph t\ge t_0$. Consider the four dimensional hamiltonian system
$$
\sist{\phi'= A(t)\phi + B(t)\psi;}{\psi' = C(t)\phi - A^*(t)\psi, \phh t\ge t_0.} \eqno (1.1)
$$
Here $\phi = (\phi_1, \phi_2), \ph \psi = (\psi_1, \psi_2)$ are the unknown continuously differentiable vector functions on $[t_0;+\infty)$. Along with the system (1.1) consider the linear system of matrix equations
$$
\sist{\Phi'= A(t)\Phi + B(t)\Psi;}{\Psi' = C(t)\Phi - A^*(t)\Psi, \phh t\ge t_0,} \eqno (1.2)
$$
Where $\Phi(t)$ and $\Psi(t)$ are the unknown continuously differentiable matrix functions of dimension $2\times 2$ on $[t_0;+\infty)$.

{\bf Definition 1.1}. {\it A solution  $(\Phi(t), \Psi(t))$ of the system (1.2) is called conjoined (or prepared, preferred) if  $\Phi^*(t)\Psi(t) = \Psi^*(t)\Phi(t), \ph t\ge t_0$.}

{\bf Definition 1.2.} {\it A solution  $(\Phi(t), \Psi(t))$ of the system (1.1) is called oscillatory if   $\det \Phi(t)$  has arbitrary large zeroes.}

{\bf Definition 1.3} {\it The system (1.1) is called oscillatory if all conjoined solutions of the system (1.2) are oscillatory, otherwise it is called non oscillatory}.

Study of the oscillatory and non oscillatory behavior of hamiltonian systems (in particular of the system (1.1)) is an important problem of qualitative theory of differential equations and many works are devoted to it (see e.g., [1 - 10] and cited works therein). For any Hermitian matrix $H$ the nonnegative (positive) definiteness of it we denote by $H \ge 0, \ph (H>0$). In the works [1 - 9] the oscillatory behavior of general hamiltonian systems is studied under the condition that the coefficient corresponding to $B(t)$ is assumed to be positive definite. In this paper we study the oscillatory and non oscillatory behavior of the system (1.1) in the direction that the assumption  $B(t) > 0, \ph t\ge t_0,$ may be destroyed.

{\bf 2. Auxiliary propositions}. Let $f(t), \ph g(t), \ph h(t), \ph h_1(t)$ be real valued continuous functions on  $[t_0;+\infty)$.  Consider the Riccati equations
$$
y' + f(t) y^2 + g(t) y + h(t) = 0, \phh t\ge t_0; \eqno (2.1)
$$
$$
y' + f(t) y^2 + g(t) y + h_1(t) = 0, \phh t\ge t_0; \eqno (2.2)
$$

{\bf Theorem 2.1}. {\it Let Eq. (2.2) has a real valued  solution $y_1(t)$  on  $[t_1;t_2) \ph (t_0\le t_1 < t_2 \le +\infty)$, and let  $f(t) \ge 0, \ph h(t) \le h_1(t), \ph t\in [t_1;t_2)$. Then for each  $y_{(0)} \ge y_1(t_0)$ Eq. (2.1)
has the solution $y_0(t)$  on  $[t_1;t_2)$   with  $y_0(t_0) = y_{(0)}$, and $y_0(t) \ge y_1(t), \ph t\in [t_1;t_2)$.}

A proof for a more general theorem is presented in [11] (see also [12]).

Denote: $I_{g,h}(\xi;t) \equiv \il{\xi}{t} \exp\biggl\{-\il{\tau}{t}g(s) d s\biggr\} h(\tau) d \tau, \ph t\ge \xi \ge t_0.$ Let $t_0 < \tau_0 \le + \infty$ and let $t_0 < t_1 < ... $ be a finite or infinite sequence such that $t_k \in [t_0;\tau_0], \ph k=1,2, ...$ We assume that  if $\{t_k\}$ is finite then  the maximum of $t_k$ is equal to $\tau_0$ and if $\{t_k\}$ is infinite then $\lim\limits_{k\to +\infty} t_k = \tau_0$.

{\bf Theorem 2.2.} {\it Let $f(t) \ge 0, \ph t\in [t_0; \tau_0), \ph t\in [t_0; \tau_0)$, and
$$
\il{t_k}{t}\exp\biggl\{\il{t_k}{\tau}\bigl[g(s) - I_{g,h}(t_k;s)\bigr]d s\biggr\} h(\tau) d \tau \le 0, \ph t\in [t_k;t_{k+1}), \ph k=0,1, ....
$$
Then for every $y_{(0)} \ge 0$ Eq. (2.1) has the solution $y_0(t)$ on $[t_0;\tau_0)$ satisfying the initial condition $y_0(t_0) = y_{(0)}$ and $y_0(t) \ge 0, \ph t\in [t_0; \tau_0)$.}

See the proof in [12].

Consider the matrix Riccati equation
$$
Z' + Z B(t) Z + A^*(t) Z + Z A(t) - C(t) = 0, \phh t\ge t_0. \eqno (2.3)
$$
The solutions $Z(t)$ of this equation existing on an interval $[t_1; t_2) (t_0 \le t_1 < t_2 \le +\infty)$ are connected with solutions $(\phi(t), \Psi(t))$ of the system (1.2) by relations (see [10]):
$$
\Phi'(t) = [A(t) + B(t) Z(t)] \Phi(t), \ph \Phi(t_1) \ne 0, \ph \Psi(t) = Z(t) \Phi(t), \ph t\in [t_1; t_2). \eqno (2.4)
$$

Let $Z_0(t)$ be a solution to Eq. (2.3) on $[t_1; t_2)$.

{\bf Definition.} {\it We will say that $[t_1; t_2)$ is the maximum existence interval for $Z_0(t)$ if $Z_0(t)$ cannot be continued to the right of $t_2$ as a solution of Eq. (2.3).}

{\bf Lemma 2.1}. {\t Let $Z_0(t)$ be a solution of Eq. (2.3) on $[t_1;t_2)$ and let $t_2 < +\infty$. Then $[t_1;t_2)$
cannot be the maximum existence interval for $Z_0(t)$ provided the function $G(t) \equiv \il{t_1}{t}tr [B(\tau) Z_0(\tau)]d\tau, \ph t\in [t_1; t_2)$, is bounded from below on $[t_1; t_2)$.}

Proof. By analogy of the proof of Lemma 2.1 from [10].

Assume $B(t) = diag \{b_1(t), b_2(t)\}, \ph t\ge t_0$. Then it
is not difficult to verify that for Hermitian unknowns $Z=\begin{pmatrix}z_{11} & z_{12}\\ \overline{z}_{12} & z_{22}\end{pmatrix}$ Eq. (2.3) is equivalent to the following nonlinear system
$$
\left\{
\begin{array}{l}
z'_{11} + b_1(t) z^2_{11} + 2 Re a_{11}(t) z_{11} + b_2(t)|z_{12}|^2 + a_{21}(t) z_{12} + \overline{a}_{21}(t) \overline{z}_{12} - c_{11}(t) = 0;\\
z'_{12} + [b_1(t) z_{11} + b_2(t) z_{22} + \overline{a}_{11}(t) + a_{22}(t)] z_{12} +  \\ \phantom{aaaaaaaaaaaaaaaaaaaaaaaaaaaaaaaaaa}+ a_{12}(t) z_{11} + a_{21}(t) z_{22} - c_{12}(t) = 0;\\
z'_{22} + b_2(t) z_{22}^2  + 2 Re a_{22}(t) z_{22} + b_1(t)|z_{12}|^2 + \overline{a}_{12}(t) z_{12} + a_{12}(t) \overline{z}_{12} - c_{22}(t) = 0,
\end{array}
\right.
\eqno (2.5)
$$
$t\ge t_0.$ If $b_2(t) \ne 0, \ph t\ge t_0,$ then it is not difficult to verify that the first equation of the system (2.5) can be rewritten  in the form
$$
z'_{11} + b_1(t) z^2_{11} + 2 Re a_{11}(t) z_{11} +  b_2(t)\left|z_{12} + \frac{\overline{a}_{21}(t)}{b_2(t)}\right|^2 - \frac{|a_{21}(t)|^2}{b_2(t)} -  c_{11}(t) = 0, \ph t\ge t_0, \eqno (2.6)
$$
and if in addition $\overline{a}_{21}(t)/ b_2(t)$ is continuously differentiable on $[t_0; +\infty)$ then by the substitution
$$
z_{12} = y - \frac{\overline{a}_{21}(t)}{b_2(t)}, \phh t \ge t_0, \eqno (2.7)
$$
in the first and second equations of the system (2.5) we get the subsystem
$$
\left\{
\begin{array}{l}
z'_{11} + b_1(t) z^2_{11} + 2 Re a_{11}(t) z_{11} +  b_2(t)|y|^2 - \frac{|a_{21}(t)|^2}{b_2(t)} -  c_{11}(t) = 0\\
 y' + [b_1(t) z_{11} + b_2(t) z_{22} + \overline{a}_{11}(t) + a_{22}(t)]y +  \bigl(a_{12}(t) - \frac{b_1(t)}{b_2(t)} \overline{a}_{21}(t)\bigr) z_{11}-\\ \phantom{aaaaaaaaaaaaaaaaaa} - \bigl(\frac{\overline{a}_{21}(t)}{b_2(t)}\bigr)' - \frac{\overline{a}_{21}(t)}{b_2(t)}\bigl(\overline{a}_{11}(t) + a_{22}(t)\bigr) - c_{12}(t) = 0, \ph t\ge t_0.
\end{array}
\right.
\eqno (2.8)
$$
Analogously if $b_1(t) \ne 0, \ph t\ge t_0,$ then the third equation of the system (2.5) can be rewritten in the form
$$
z'_{22} + b_2(t) z^2_{22} + 2 Re a_{22}(t) z_{22} +  b_1(t)\left|z_{12} + \frac{a_{12}(t)}{b_1(t)}\right|^2 - \frac{|a_{12}(t)|^2}{b_1(t)} -  c_{22}(t) = 0, \ph t\ge t_0, \eqno (2.9)
$$
and if in addition $a_{12}(t)/ b_1(t)$ is continuously differentiable on $[t_0; +\infty)$ then by the substitution
$$
z_{12} = v - \frac{a_{12}(t)}{b_1(t)}, \phh t \ge t_0, \eqno (2.10)
$$
in the second and third equations of the system (2.5) we obtain the subsystem
$$
\left\{
\begin{array}{l}
z'_{22} + b_2(t) z^2_{22} + 2 Re a_{22}(t) z_{22} +  b_1(t)|v|^2 - \frac{|a_{12}(t)|^2}{b_1(t)} -  c_{22}(t) = 0\\
 y' + [b_1(t) z_{11} + b_2(t) z_{22} + \overline{a}_{11}(t) + a_{22}(t)]v +  \bigl(\overline{a}_{21}(t) - \frac{b_2(t)}{b_1(t)} a_{12}(t)\bigr) z_{22}-\\ \phantom{aaaaaaaaaaaaaaaaa} - \bigl(\frac{a_{12}(t)}{b_1(t)}\bigr)' - \frac{a_{12}(t)}{b_1(t)}\bigl(\overline{a}_{11}(t) + a_{22}(t)\bigr) - c_{12}(t) = 0, \ph t\ge t_0.
\end{array}
\right.
\eqno (2.11)
$$
If $(z_{11}(t), y(t))$ is a solution of the subsystem (2.8) on $[t_0;t_1) (t_0 < t_1 \le + \infty)$ with $y(t_0) = 0$ and $(z_{22}(t), v(t))$ is a solution of the subsystem (2.11) on $[t_0;t_1)$ with $v(t_0)=0$   then by Cauchi formula from the second equation of the subsystem (2.8) and from the second equation of the subsystem (2.11) we have respectively:
$$
y(t) = - \exp\biggl\{-\il{t_0}{t}b_1(\tau)z_{11}(\tau)d\tau\biggr\}\il{t_0}{t}\biggl[\exp\biggl\{\il{t_0}{\tau}b_1(s) z_{11}(s)d s\biggr\}\biggr]'\biggl(\frac{a_{12}(\tau)}{b_1(\tau)} - \frac{\overline{a}_{21}(\tau)}{b_2(\tau)}\biggr)\times
$$
$$
\times\exp\biggl\{-\il{\tau}{t}\bigl(b_2(s) z_{22}(s) + \overline{a}_{11}(s) + a_{22}(s)\bigr)ds\biggr\}d\tau +
$$
$$
+\il{t_0}{t}\exp\biggl\{-\il{\tau}{t}\bigl(b_1(s)z_{11}(s) + b_2(s) z_{22}(s) + \overline{a}_{11}(s) + a_{22}(s)\bigr)ds\biggr\}\biggl[\biggl(\frac{\overline{a}_{21}(\tau)}{b_2(\tau)}\biggr)' +\phantom{aaaaaaaaa}
$$
$$
\phantom{aaaaaaaaaaaaaaaaaaaaaaaaaaaaaaaaaaaaaa}+\frac{\overline{a}_{21}(\tau)}{b_2(\tau)}\bigl(\overline{a}_{11}(\tau) + a_{22}(\tau)\bigr) + c_{12}(\tau)\biggr]d\tau,
$$
$$
v(t) = - \exp\biggl\{-\il{t_0}{t}b_2(\tau)z_{22}(\tau)d\tau\biggr\}\il{t_0}{t}\biggl[\exp\biggl\{\il{t_0}{\tau}b_2(s) z_{22}(s)d s\biggr\}\biggr]'\biggl(\frac{\overline{a}_{21}(\tau)}{b_2(\tau)} - \frac{a_{12}(\tau)}{b_1(\tau)} \biggr)\times
$$
$$
\times\exp\biggl\{-\il{\tau}{t}\bigl(b_1(s) z_{11}(s) + \overline{a}_{11}(s) + a_{22}(s)\bigr)ds\biggr\}d\tau +
$$
$$
+\il{t_0}{t}\exp\biggl\{-\il{\tau}{t}\bigl(b_1(s)z_{11}(s) + b_2(s) z_{22}(s) + \overline{a}_{11}(s) + a_{22}(s)\bigr)ds\biggr\}\biggl[\biggl(\frac{a_{12}(\tau)}{b_1(\tau)}\biggr)' +\phantom{aaaaaaaaa}
$$
$$
\phantom{aaaaaaaaaaaaaaaaaaaaaaaaaaa}+\frac{a_{12}(\tau)}{b_1(\tau)}\bigl(\overline{a}_{11}(\tau) + a_{22}(\tau)\bigr) + c_{12}(\tau)\biggr]d\tau, \ph t\in [t_0;t_1).
$$
From here it is easy to derive

{\bf Lemma 2.2}. {\it Let $b_j(t) > 0, \ph j=1,2,$ the functions $a_{12}(t)/b_1(t), \ph \overline{a}_{21}(t)/b_2(t)$ be continuously differentiable on $[t_0; t_1) (t_0 < t_1 < + \infty))$ and let  $(z_{11}(t), y(t))$ and $(z_{22}(t), v(t))$ be solutions of the subsystems (2.8) and (2.11) respectively on  $[t_0; t_1)$ such that $z_{jj}(t) \ge~ 0, \linebreak t\in ~ [t_0;t_1), \ph j=1,2, \ph y(t_0) = v(t_0) = 0$. Then
$$
|y(t)| \le \mathfrak{M}(t) + \il{t_0}{t}\biggl|\exp\biggl\{-\il{\tau}{t}\bigl(\overline{a}_{11}(s) + a_{22}(s)\bigr)ds\biggr\}\biggl[\biggl(\frac{\overline{a}_{21}(\tau)}{b_2(\tau)}\biggr)'+
$$
$$
\phantom{aaaaaaaaaaaaaaaaaaaaaaaaaaaaaaaaaaaaaaa}+\frac{\overline{a}_{21}(\tau)}{b_2(\tau)}\bigl(\overline{a}_{11}(\tau) + a_{22}(\tau)\bigr)+ + c_{12}(\tau)\biggr]\biggr|d\tau,
$$

$$
|v(t)| \le \mathfrak{M}(t) + \il{t_0}{t}\biggl|\exp\biggl\{-\il{\tau}{t}\bigl(\overline{a}_{11}(s) + a_{22}(s)\bigr)ds\biggr\}\biggl[\biggl(\frac{a_{12}(\tau)}{b_1(\tau)}\biggr)'+
$$
$$
\phantom{aaaaaaaaaaaaaaaaaaaaaaaaaaaa}+\frac{a_{12}(\tau)}{b_1(\tau)}\bigl(\overline{a}_{11}(\tau) + a_{22}(\tau)\bigr)+ + c_{12}(\tau)\biggr]\biggr|d\tau, \ph t\in [t_0;t_1),
$$
where
$$
\mathfrak{M}(t)\equiv \max\limits_{\tau\in [t_0; t]}\biggl|\exp\biggl\{-\il{\tau}{t}\bigl(\overline{a}_{11}(s) + a_{22}(s)\bigr)ds\biggr\}\biggl(\frac{a_{12}(\tau)}{b_1(\tau)} - \frac{\overline{a}_{21}(\tau)}{b_2(\tau)}\biggr)\biggr|, \ph t\ge t_0.
$$
}
\phantom{aaaaaaaaaaaaaaaaaaaaaaaaaaaaaaaaaaaaaaaaaaaaaaaaaaaaaaaaaaaaaaaaaaaaaaaaaa}$\Box$

{\bf Lemma 2.3.} {\it For any two square matrices $M_1\equiv (m_{ij}^1)_{ij=1}^n, \ph M_2\equiv (m_{ij}^2)_{ij=1}^n$ the equality
$$
tr (M_1 M_2) = tr (M_2 M_1)
$$
is valid.}

Proof. We have: $tr (M_1 M_2) = \sum\limits_{j=1}^n(\sum\limits_{k=1}^n m_{jk}^1 m_{kj}^2) = \sum\limits_{k=1}^n(\sum\limits_{j=1}^n m_{jk}^1 m_{kj}^2) = \sum\limits_{k=1}^n(\sum\limits_{j=1}^n m_{kj}^2 m_{jk}^1) = tr (M_2 M_1).$ The lemma is proved.

{\bf 3. Main results}. Let $f_{jk}(t), \ph j,k =1,2, \ph t\ge t_0,$ be real valued continuous functions on $[t_0; +\infty)$. Consider the linear system of equations
$$
\sist{\phi_1' = f_{11}(t) \phi_1 + f_{12}(t) \psi_1;}{\psi_1' = f_{21}(t) \phi_1 + f_{22}(t) \psi_1, \ph t\ge t_0,} \eqno (3.1)
$$
and the Riccati equation
$$
y' + f_{12}(t) y^2 + [f_{11}(t) - f_{22}(t)] y - f_{12}(t) = 0, \phh t\ge t_0. \eqno (3.2)
$$
All solutions $y(t)$ of the last equation, existing on some interval $[t_1; t_2)\hskip 2pt (t_0 \le t_1 < t_2 \le + \infty)$ are connected with solutions $(\phi_1(t), \psi_1(t))$ of the system (3.1) by relations (see [13]):
$$
\phi_1(t) = \phi_1(t_1)\exp\biggl\{\il{t_1}{t}\bigl[f_{12}(\tau) y(\tau) + f_{11}(\tau)\bigr]d\tau\biggr\}, \ph \phi_1(t_1)\ne 0, \ph \psi_1(t) = y(t) \phi_1(t),    \eqno (3.3)
$$
$t\in [t_1; t_2).$

{\bf Definition 3.1.} {\it The system (3.1) is called oscillatory if for its every solution \linebreak $(\phi_1(t), \psi_1(t))$ the function $\phi_1(t)$ has arbitrary large zeroes.}

{\bf Remark 3.1.} {\it Some explicit oscillatory criteria for the system (3.1) are proved in [10] amd [14]}.

{\bf 3.1. The case when $B(t)$ is a diagonal matrix}. In this subsection we will assume that $B(t) = diag\{b_1(t), b_2(t)\}$.
Denote:
$$
\chi_j(t) \equiv \sist{c_{jj}(t) \ph if b_{3-j}(t) = 0;}{c_{jj}(t) + \frac{|a_{3-j,j}(t)|^2}{b_{3-j}(t)}, \ph if  b_{3-j}(t) \ne 0,} \phh t\ge t_0, \ph j=1,2.
$$

{\bf Theorem 3.1.} {\it Assume $b_j(t) \ge 0, \ph t\ge t_0,$ and if $b_j(t) = 0$ then $a_{3-j,j}(t) = 0, \ph j=1,2, \ph t\ge t_0.$
Under these restrictions the system (1.1) is oscillatory provided one of the systems
$$
\sist{\phi_1'= 2 Re (a_{jj}(t)) \phi_1 + b_j(t) \psi_1;}{\psi_1' = - \chi_j(t) \phi_1, \phh t\ge t_0,} \eqno (3.4_j)
$$
j=1,2, is oscillatory.
}

Proof. Suppose the system (1.1) is not oscillatory. Then for some conjoined solution $(\Phi(t), \Psi(t))$ of the system (1.2) there exists $t_1 \ge t_0$ such that $det \Phi(t) \ne 0, \ph t\ge t_1.$ Due to (2.4) from here it follows that $Z(t)\equiv \Psi(t) \Phi^{-1}(t), \ph  t\ge t_1,$ is a Hermitian solution to Eq. (2.3) on $[t_1; +\infty)$.
Let $Z(t) = \begin{pmatrix}z_{11}(t) & z_{12}(t)\\ \overline{z}_{12}(t) & z_{22}(t)\end{pmatrix}, \ph t\ge t_1.$
Consider the Riccati equations
$$
y' + b_1(t) y^2 + 2 (Re a_{11}(t)) y + b_2(t)|z_{12}(t)|^2 + a_{21}(t) z_{12}(t) + \overline{a}_{21}(t) \overline{z}_{12}(t) - c_{11}(t) = 0, \eqno (3.5)
$$
$$
y' + b_2(t) y^2 + 2 (Re a_{22}(t)) y + b_1(t)|z_{12}(t)|^2 + \overline{a}_{12}(t) z_{12}(t) + a_{12}(t) \overline{z}_{12}(t) - c_{22}(t) = 0, \eqno (3.6)
$$
$$
y' + b_j(t) y^2 + 2 (Re a_{jj}(t) y + \chi_j(t) = 0, \eqno (3.7_j)
$$
$j=1,2, \ph t\ge t_1.$ By (2.6) and (2.9) from the conditions of the theorem it follows that
$$
\chi_1(t) \le b_2(t)|z_{12}(t)|^2 + a_{21}(t) z_{12}(t) + \overline{a}_{21}(t) \overline{z}_{12}(t) - c_{11}(t), \ph t\ge t_1,
$$
$$
\chi_2(t) \le  b_1(t)|z_{12}(t)|^2 + \overline{a}_{12}(t) z_{12}(t) + a_{12}(t) \overline{z}_{12}(t) - c_{22}(t), \ph t\ge t_1.
$$
Using Theorem 2.1 to the pairs (3.5), $(3.7_1)$ and (3.6), $(3.7_2)$ of equations from here we conclude that the equations  $(3.7_j), \ph j=1,2,$ have solutions on $[t_1; +\infty)$. By (3.1) - (3.3) from here it follows that the systems  $(3.4_j), \ph j=1,2,$ are not oscillatory which contradicts the condition of the theorem. The obtained contradiction completes the proof of the theorem.

Denote:
$
I_j(\xi;t) \equiv \il{\xi}{t}\exp\biggl\{-\il{\tau}{t}2(Re a_{jj}(s))d s \biggr\}\chi_j(\tau) d\tau, \ph t\ge \xi \ge t_0, \ph j=1,2.
$

{\bf Theorem 3.2.} {\it Assume $b_1(t) \ge 0 (\le 0), \ph b_2(t) \le 0 (\ge 0)$ and if $b_j(t) = 0$ then $a_{j, 3-j}(t) = 0, \ph j=1,2, \ph t\ge t_0$; there exist infinitely large sequences $\xi_{j,0} = t_) < \xi_{j,1}  < ... < \xi_{j,m} , ..., \ph j=1,2, $ such that
$$
1_j) \ph (-1)^j \il{\xi_{j,m}}{t} \exp\biggl\{\il{\xi_{j,m}}{\tau}\biggl[2 Re  a_{jj}(s) - (-1)^j I_j(xi_{j,m},s)\biggr] d s\biggr\} \chi_j(\tau) d\tau \ge 0 \ph (\le 0),
$$
$\ph t\in [\xi_{j,m}; \xi_{j,m +1}), \ph    m=1,2,3, ....., \ph j=1,2$. Then the system (1.1) is non oscillatory.}

Proof. Let us prove the theorem only in the case when $b_1(t) \ge 0, \ph b_2(t) \le 0, \ph t\ge~ t_0$. The case $b_1(t) \le 0, \ph b_2(t) \ge 0, \ph t\ge t_0$, can be proved by analogy. Let $(\Phi(t), \Psi(t))$ be a conjoined solution of the system (1.2) with $\Phi(t_0) = \begin{pmatrix}1 & 0\\ 0 & -1\end{pmatrix}$ and let $[t_0; T)$ be the maximum interval such that $det \Phi(t) \ne 0, \ph t\in [t_0; T)$. Then by (2.4) the matrix function $Z(t) \equiv \Psi(t) \phi^{-1}(t), \ph t\in [t_0; T)$, is a Hermitian solution to Eq. (2.3) on $[t_0; T)$. By (2.5), (2.7), (2.8), (2.10), (2.11) from here it follows that the subsystems (2.8) and (2.11) have solutions $(z_{11}(t), y(t))$ and $(z_{22}(t), v(t))$ respectively on $[t_0; T)$ with $z_{11}(t_0) = 1, \linebreak z_{22}(t_0) =-1$. Show that
$$
z_{11}(t) \ge 0, \phh t\in [t_0; T). \eqno (3.8)
$$
Consider the Riccati equations
$$
z' + b_1(t) z^2 + 2 (Re a_{11}(t)) z + b_2(t)|y(t)|^2 + \chi_1(t) = 0, \phh t\in [t_0;T), \eqno (3.9)
$$
$$
z' + b_1(t) z^2 + 2 (Re a_{11}(t)) z  +  \chi_1(t) = 0, \phh t\in [t_0;T), \eqno (3.10)
$$
By Theorem 2.2 from the conditions of the theorem it follows that the last equation has a nonnegative solution on $[t_0; T).$ Then using Theorem 2.1 to the pair of equations (3.9), (3.10) on the basis of the conditions of the theorem we conclude that Eq. (3.9) has a nonnegative solution $z_0(t)$  on $[t_0; T)$ with $z_0(t_0) = 0$. Then since $z_{11}(t)$ is a solution to Eq. (3.9) on $[t_0;T)$ and $z_{11}(t_0) =1$ we have (3.8). Show that
$$
z_{22}(t) \le 0, \phh t\in [t_0; T). \eqno (3.11)
$$
Consider the Riccati equations
$$
z' - b_2(t) z^2 + 2 (Re a_{22}(t)) z   - \chi_2(t) = 0, \phh t\in [t_0;T), \eqno (3.12)
$$
$$
z' - b_2(t) z^2 + 2 (Re a_{22}(t)) z  - b_1(t)|v(t)|^2  - \chi_2(t) = 0, \phh t\in [t_0;T). \eqno (3.13)
$$
By Theorem 2.2 from the conditions of the theorem it follows that Eq. (3.12) has a nonnegative solution $z_1(t)$ on $[t_0; T)$ with $z_1(t_0) = 0$. Then using Theorem 2.1 to the pair of equations (3.12) and (3.13) we derive that Eq. (3.13) has a nonnegative solution $z_2(t)$ on $[t_0; T)$ whit $z_2(t_0) = 0$. Hence since obviously $-z_{22}(t)$ is a solution of Eq. (3.13) on $[t_0; T)$ and $-z_{11}(t_0) =1$ we have (3.11). Since $b_1(t) \ge 0, \ph b_2(t) \le 0, \ph t\in [t_0;T)$ from (3.8) and (3.11) it follows:
$$
\il{t_0}{t}\Bigl[b_1(\tau) z_{11}(\tau) + b_2(\tau) z_{22}(\tau)\Bigr]d\tau \ge 0, \phh t\in [t_0; T). \eqno (3.14)
$$
To complete the proof of the theorem it remains to show that $T = + \infty$. Suppose $T < + \infty$. Then by virtue of Lemma 2.1 from (3.14) it follows that $[t_0; T)$ is not the maximum existence interval for $Z(t)$. By (2.4) from here it follows that $det \Phi(t) \ne 0, \ph t\in [t_0; T_1),$ for some $T_1 > T$. We have obtained a contradiction which completes the proof of the theorem.

{\bf Remark 3.2.} {\it The conditions $1_j), \ph j=1,2,$ are satisfied if in particular $(-1)^j\chi_j(t) \ge~ 0 \linebreak (\le 0), \ph t \ge t_0$.}

Denote:
$$
\chi_3(t) \equiv b_2(t)\biggl[\mathfrak{M}(t) + \il{t_0}{t}\biggl|\exp\biggl\{- \il{\tau}{t}\bigl[\overline{a}_{11}(s) + a_{22}(s)\bigr]d s\biggr\}\times \phantom{aaaaaaaaaaaaaaaaaaaaaaaaaaaaaaaaaaa}
$$
$$
\phantom{aa}\times \biggl[\biggl(\frac{\overline{a}_{21}(t)}{b_2(t)}\biggl)' + \frac{\overline{a}_{21}(\tau)}{b_2(\tau)}(\overline{a}_{11}(\tau) + a_{22}(\tau)\biggr) + c_{12}(\tau)\biggr]\biggr|d\tau\biggr]^2 - \frac{|a_{21}(t)|^2}{b_2(t)} - c_{11}(t),
$$

$$
\chi_4(t) \equiv b_1(t)\biggl[\mathfrak{M}(t) + \il{t_0}{t}\biggl|\exp\biggl\{- \il{\tau}{t}\bigl[\overline{a}_{11}(s) + a_{22}(s)\bigr]d s\biggr\}\times \phantom{aaaaaaaaaaaaaaaaaaaaaaaaaaaaaaaaaaa}
$$
$$
\phantom{aa}\times \biggl[\biggl(\frac{a_{12}(t)}{b_1(t)}\biggl)' + \frac{a_{12}(\tau)}{b_1(\tau)}(\overline{a}_{11}(\tau) + a_{22}(\tau)\biggr) + c_{12}(\tau)\biggr]\biggr|d\tau\biggr]^2 - \frac{|a_{12}(t)|^2}{b_1(t)} - c_{22}(t),
$$
$$
I_{j+2}(\xi;t) \equiv \il{\xi}{t}\exp\biggl\{-\il{\tau}{t}2(Re a_{jj}(s))d s \biggr\}\chi_{j+2}(\tau) d\tau, \ph  t\ge \xi \ge t_0, \ph j=1,2.
$$

{\bf Theorem 3.3.} {\it Let the following conditions be satisfied

\noindent
1) $b_j(t)> 0, \ph t\ge t_0, \ph j=1,2;$

\noindent
2) the functions $a_{12}(t)/b_1(t)$ and $\overline{a}_{21}(t)/b_2(t)$ are continuously differentiable on $[t_0;+\infty)$;

\noindent
3) there exist infinitely large sequences $\xi_{j,0} = t_0 < \xi_{j,1}  < ... < \xi_{j,m} , ..., \ph j=1,2, $ such that
$$
\il{\xi_{j,m}}{t} \exp\biggl\{\il{\xi_{j,m}}{\tau}\biggl[2 Re  a_{jj}(s) -  I_{j+2}(\xi_{j,m},s)\biggr] d s\biggr\} \chi_{j+2}(\tau) d\tau \le 0, \ph t\in [\xi_{j,m}; \xi_{j,m +1}),
$$
$m=1,2,3, ....., \ph j=1,2$. Then the system (1.1) is non oscillatory.}

Proof. Let $Z(t) \equiv \begin{pmatrix}z_{11}(t) & z_{12}(t)\\\overline{z}_{12}(t) & z_{22}(t)\end{pmatrix}$ be the Hermitian solution of Eq. (2.3) on $[t_0; T)$ satisfying the initial condition $Z(t_0) = \begin{pmatrix}1 & 0\\0 & 1\end{pmatrix}$, where $[t_0;T)$ is the maximum existence interval for $Z(t)$. Due to (2.4) to prove the theorem it is enough to show that
$$
T= + \infty. \eqno (3.15)
$$
By (2.5), (2.7), (2.8), (2.10), (2.11) from the conditions 1) and 2) it follows that \linebreak $(z_{11}(t), z_{12}(t) + \overline{a}_{21}(t)/ b_2(t))$ and $(z_{22}(t), z_{12}(t) + a_{12}(t)/ b_1(t))$ are solutions of the subsystems (2.8) and (2.11) respectively on $[t_0; T)$. Show that
$$
z_{jj}(t) > 0, \phh t\in [t_0; T). \eqno (3.16)
$$
Suppose it is not so. Then there exists $T_1 \in (t_0; T)$ such that
$$
z_{11}(t)z_{22}(t) > 0, \ph t\in [t_0; T_1), \ph z_{11}(T_1)z_{22}(T_1) = 0. \eqno (3.17)
$$
Without loss of generality we may take that $a_{12}(t_0) = a_{21}(t_0) = 0$. Then by virtue of Lemma 2.2 from (3.17) it follows that
$$
\biggl|z_{12}(t) + \frac{\overline{a}_{21}(t)}{b_2(t)}\biggr| \le \mathfrak{M}(t) + \il{t_0}{t}\biggl|\exp\biggl\{ - \il{\tau}{t}\bigl(\overline{a}_{11}(s) + a_{22}(s)\biggr)d s\biggr\}\biggl[\biggl(\frac{\overline{a}_{21}(\tau)}{b_2(\tau)}\biggr)' + \phantom{aaaaaaaaaaaaaaaaaaaaaaaaaaaaaa}
$$
$$
\phantom{aaaaaaaaaaaaaaaaaaaaaaaaaaaaaaaa} + \frac{\overline{a}_{21}(\tau)}{b_2(\tau)}\bigl(\overline{a}_{11}(\tau) + a_{22}(\tau)\bigr) - c_{12}(\tau)\biggr]\biggr|d\tau,
$$
$$
\biggl|z_{12}(t) + \frac{a_{12}(t)}{b_1(t)}\biggr| \le \mathfrak{M}(t) + \il{t_0}{t}\biggl|\exp\biggl\{ - \il{\tau}{t}\bigl(\overline{a}_{11}(s) + a_{22}(s)\biggr)d s\biggr\}\biggl[\biggl(\frac{a_{12}(\tau)}{b_1(\tau)}\biggr)' + \phantom{aaaaaaaaaaaaaaaaaaaaaaaaaaaaaaaaa}
$$
$$
\phantom{aaaaaaaaaaaaaaaaaaaaa} + \frac{a_{12}(\tau)}{b_1(\tau)}\bigl(\overline{a}_{11}(\tau) + a_{22}(\tau)\bigr) - c_{12}(\tau)\biggr]\biggr|d\tau, \ph t\in [t_0;T_1).
$$
Hence
$$
b_2(t)\biggl|z_{12}(t) + \frac{\overline{a}_{21}(t)}{b_2(t)}\biggr| - \frac{|a_{21}(t)|^2}{b_2(t)} - c_{11}(t) \le \chi_3(t), \phantom{aaaaaaaaaaaaaaaaaaaaaaaaaaaaaaaaaaaaaaaaaaaaaaaaaa}
$$
$$
\phantom{aaaaaaaaaaaaaa}b_1(t)\biggl|z_{12}(t) + \frac{a_{12}(t)}{b_1(t)}\biggr|^2 -  \frac{|a_{12}(t)|^2}{b_2(t)} - c_{22}(t) \le \chi_4(t), \phh t\in [t_0;T_1),
$$
By virtue of Theorem 2.1 and Theorem 2.2 from here and from the condition 3) it follows that the Riccati equations
$$
z' + b_1(t) z^2 + 2(Re a_{11}(t)) z + b_2(t)\biggl|z_{12}(t) + \frac{\overline{a}_{21}(t)}{b_2(t)}\biggr| - \frac{|a_{21}(t)|^2}{b_2(t)} - c_{11}(t) = 0,  \eqno (3.18)
$$
$$
z' + b_2(t) z^2 + 2(Re a_{22}(t)) z + b_1(t)\biggl|z_{12}(t) + \frac{a_{12}(t)}{b_1(t)}\biggr|^2 -  \frac{|a_{12}(t)|^2}{b_2(t)} - c_{22}(t) = 0,  \eqno (3.19)
$$
$ t\in [t_0; T_1),$ have nonnegative solutions $z_1(t)$ and $z_2(t)$  respectively on $[t_0; T_1)$ with $z_1(t_0) = z_2(t_0) = 0$. Obviously $z_{11}(t)$ and $z_{22}(t)$ are solutions of Eq. (3.18) and (3.19) respectively on $[t_0; T_1]$. Therefore since $z_{jj}(t_0) = 1 > z_j(t_0) = 0, \ph j=1,2$ due to uniqueness theorem $z_{jj}(t) > 0, \ph t\in [t_0;T_1], \ph j=1,2,$ which contradicts (3.17). The obtained contradiction proves (3.16). From (3.16) and 1) it follows that
$$
\il{t_0}{t}\bigl[b_1(\tau) z_{11}(\tau) + b_2(\tau) z_{22}(\tau)\bigr] d \tau \ge 0, \phh t\in [t_0; T). \eqno (3.20)
$$
Suppose $T < + \infty$. Then by Lemma 2.1 from (3.20) it follows that $[t_0; T)$ is not the maximum existence interval for $Z(t)$ which contradicts our assumption. The obtained contradiction proves (3.15). The theorem is proved.

{\bf Remark 3.3.} {\it The conditions 3) of Theorem 3.3 are satisfied if in particular $\chi_j(t) \le ~0, \linebreak t\ge t_0, \ph j=1,2.$}

{\bf 3.2. The case when $B(t)$ is nonnegative definite}.
In this subsection we will assume that $B(t)$ is nonnegative definite and  $\sqrt{B(t)}$ is continuously differentiable on $[t_0;+ \infty)$. Consider the matrix equation
$$
\sqrt{B(t)} X [A(t) \sqrt{B(t)} - \sqrt{B(t)}'] = A(t) \sqrt{B(t)} - \sqrt{B(t)}', \phh t\ge t_0. \eqno (3.21)
$$
Obviously this equation has always a solution on $[a;b] (\subset [t_0; + \infty))$ when $B(t) > 0, \linebreak t\in [a;b] \ph (X(t) = B^{-1}(t), \ph t\in [a;b])$. It may have also a solution on $[a;b]$ in some cases when $B(t) \ge 0, \ph t\in [a;b]$ (e.g., $A(t) = \begin{pmatrix} a_1(t) & a_2(t)\\ 0 & 0 \end{pmatrix},  \ph  B(t) = \begin{pmatrix} b_1(t) &  0 \\ 0 & 0 \end{pmatrix}, \ph b_1(t) >~ 0, \linebreak t\in [a;b]).$  In this subsection we also will assume that Eq. (3.21) has always a solution on $[t_0; + \infty)$.  Let $F(t)$ be a solution of Eq. (3.21) on $[t_0; + \infty)$.  Denote:
$$P(t) \equiv F(t) [A(t)\sqrt{B(t)} - \sqrt{B(t)}'] = (p_{jk}(t))_{j,k =1}^2,  \eqno (3.22)$$
$\ph Q(t) \equiv \sqrt{B(t)} C(t) \sqrt{B(t)}, \ph  (q_{jk}(t))_{j,k =1}^2, \ph \widetilde{\chi}_j(t) \equiv q_{jj}(t) + |p_{3 -j, j}(t)|^2, \ph j=~1,2, \ph t\ge t_0$.

{\bf Corollary 3.1}. {\it The system (1.1) is oscillatory provided one of the equations
$$
\phi_1'' + 2 [Re p_{jj}(t)] \phi_1' + \widetilde{\chi}_j(t) \phi_1 = 0, \ph j =1,2, \ph t\ge t_0. \eqno (3.23_j)
$$
is oscillatory.}

Proof. Multiply Eq. (2.3) at left and at right by $\sqrt{B(t)}$. Taking into account the equality $(\sqrt{B(t)}Z \sqrt{B(t)})' = \sqrt{B(t)}Z' \sqrt{B(t)} + \sqrt{B(t)}'Z \sqrt{B(t)} + \sqrt{B(t)}Z \sqrt{B(t)}' \ph t\ge t_0,$
we obtain
$$
V' + V^2 + P^*(t) V + V P(t) - Q(t) = 0, \ph t\ge t_0, \eqno (3.24)
$$
where $V \equiv \sqrt{B(t)}Z \sqrt{B(t)}$. To this equation corresponds the following matrix hamiltonian system
$$
\sist{\Phi'= P(t)\Phi + \Psi;}{\Psi' = Q(t)\Phi - P^*(t)\Psi, \phh t\ge t_0.} \eqno (3.25)
$$
Suppose the system (1.1) is not oscillatory. Then by (2.4)  Eq. (2.3) has a Hermitian solution $Z(t)$ on $[t_1; + \infty)$ for some $t_1 \ge t_0$. Therefore $V(t) \equiv \sqrt{B(t)} Z(t) \sqrt{B(t)}, \ph t\ge t_1,$ is a hermitian solution of Eq. (3.24) on $[t_1; + \infty)$ and hence the system (3.25) has a conjoined solution $(\Phi (t), \Psi(t))$ such that $det \Phi(t) \ne 0, \ph t\ge t_.$ It means that the hamiltonian system
$$
\sist{\phi'= P(t)\phi + \psi;}{\psi' = Q(t)\phi - P^*(t)\psi, \phh t\ge t_0,}
$$
is not oscillatory. By Theorem 3.1 from here it follows that the scalar systems
$$
\sist{\phi_1'= 2 Re p_{jj}(t)\phi_1  + \psi_1;}{\psi_1' = - \widetilde{\chi}_j(t) \phi_1, \phh t\ge t_0,}
$$
$j = 1,2,$ are not oscillatory. Therefore the corresponding equations $(3.23_j), \ph j=1,2,$ are not oscillatory, which contradicts the conditions of the corollary. The corollary is proved.

Denote:
$$
\widetilde{\mathfrak{M}}(t)\equiv \max\limits_{\tau\in [t_0;t]}\biggl|\exp\biggl\{-\il{\tau}{t}\bigl(\overline{p}_{11}(s) + p_{22}(s)\bigr)ds\biggr\}(p_{12}(\tau) - \overline{p}_{21}(\tau))\biggr|;
$$
$$
\widetilde{\chi}_3(t) \equiv \biggl[\widetilde{\mathfrak{M}}(t) + \il{t_0}{t}\biggl|\exp\biggl\{- \il{\tau}{t}[\overline{p}_{11}(s) + p_{22}(s)]d s\biggr\}\times \phantom{aaaaaaaaaaaaaaaaaaaaaaaaaaaaaaaaaaa}
$$
$$
\phantom{aaaaaaaaaaaaaaa}\times \Bigl[\overline{p}_{21} \hskip 0.1pt'(t) + \overline{p}_{21}(\tau)(\overline{p}_{11}(\tau) + p_{22}(\tau)) + q_{12}(\tau)\Bigr]\biggr|d\tau\biggr]^2 - |p_{21}(t)|^2 - q_{11}(t);
$$
$$
\widetilde{\chi}_4(t) \equiv \biggl[\widetilde{\mathfrak{M}}(t) + \il{t_0}{t}\biggl|\exp\biggl\{- \il{\tau}{t}\bigl[\overline{p}_{11}(s) + p_{22}(s)\bigr]d s\biggr\}\times \phantom{aaaaaaaaaaaaaaaaaaaaaaaaaaaaaaaaaaa}
$$
$$
\phantom{aaaaaa}\times \Bigl[p_{12}\hskip0.1pt'(t) + p_{12}(\tau)(\overline{p}_{11}(\tau) + p_{22}(\tau)) + q_{12}(\tau)\Bigr]\biggr|d\tau\biggr]^2 - |p_{12}(t)|^2 - q_{22}(t), \ph t\ge t_0;
$$
$$
\widetilde{I}_{j+2}(\xi,t) \equiv \il{\xi}{t}\exp\biggl\{-\il{\tau}{t}2(Re\hskip 2pt p_{jj}(s))d s \biggr\}\widetilde{\chi}_{j+2}(\tau) d\tau, \ph  t\ge \xi \ge t_0,  \ph j=1,2.
$$

{\bf Theorem 3.4.} {\it Let the following conditions be satisfied:

\noindent
$1') B(t) \ge 0, \ph t\ge t_0;$

\noindent
$2')$ Eq. (3.21) has a solution $F(t)$ on $[t_0; + \infty)$

\noindent
$3')$ the functions $p_{12}(t)$ and $p_{21}(t)$, defined by (3.22) are continuously differentiable on $[t_0; + \infty)$;

\noindent
4') there exist infinitely large sequences $\xi_{j,0} = t_0 < \xi_{j,1} < ... < \xi_{j,m}, ...$ such that
$$
\il{\xi_{j,m}}{t} \exp\biggl\{\il{\xi_{j,m}}{\tau}\biggl[2 Re \hskip 2pt  a_{jj}(s) -  \widetilde{I}_{j+2}(\xi_{j,m},s)\biggr] d s\biggr\} \widetilde{\chi}_{j+2}(\tau) d\tau \le 0, \ph t\in [\xi_{j,m}; \xi_{j,m +1}),
$$
$ \ph    m=1,2,3, ....., \ph j=1,2$. Then the system (1.1) is non oscillatory.}

Proof. Let $Z(t) \equiv \begin{pmatrix} z_{11}(t) & z_{12}(t)\\ \overline{z}_{12}(t) & z_{22}(t)\end{pmatrix}$ be the Hermitian solution of Eq. (2.3) satisfying the initial condition $Z(t_0) = \begin{pmatrix} 1 & 0\\ 0 & 1\end{pmatrix}$, and let $[t_0;T)$ be the maximum existence interval for $Z(t)$. Then $V(t) \equiv \sqrt{B(t)}Z(t) \sqrt{B(t)}$ is a soluyion of Eq. (3.24) on $[t_0; T)$. Without loss of generality we may assume that $B(t_0) = \begin{pmatrix} 1 & 0\\ 0 & 1\end{pmatrix}$. Then $V(t_0) = \begin{pmatrix} 1 & 0\\ 0 & 1\end{pmatrix}$, and by analogy of the proof of Theorem 3.3 we can show that from the conditions of the theorem it follows that
$$
\il{t_0}{t} tr V(\tau) d \tau \ge 0, \phh t\in [t_0;T). \eqno (3.26)
$$
By virtue of Lemma 2.3 we have: $tr V(t) = tr [B(t) Z(t)], \ph t\in [t_0;T)$. From here and from (3.26) it follows:
$$
\il{t_0}{t} tr [B(\tau) Z(\tau)] d \tau \ge 0, \phh t\in [t_0;T). \eqno (3.27)
$$
To complete the proof of the theorem it remains to show that $T = + \infty$. Suppose $T < + \infty$. Then by virtue of Lemma 2.2 from (3.27) it follows that $[t_0;T)$ is not the maximum existence interval for $Z(t)$ which contradicts our assumption. The obtained contradiction shows that $T = + \infty$. The theorem is proved.

Example 3.1. Consider the second order vector equation
$$
\phi'' + K(t)\phi = 0, \phh t\ge t_0, \eqno (3.28)
$$
where $K(t) \equiv \begin{pmatrix} \mu(t) & 10 i\\ -10 i & - t^2\end{pmatrix}, \ph \mu(t) \equiv p_1 \sin (\lambda_1 t + \theta_1) + p_2 \sin (\lambda_2 t + \theta_2), \ph t\ge t_0, \ph  p_j, \ph \lambda_j\ne~ 0, \linebreak \theta_j, \ph j=1,2,$ are some real constants such that $\lambda_1$ and $\lambda_2$ are rational independent. This equation is equivalent to the system (1.1) with $A(t)\equiv~ 0, \ph B(t) \equiv \begin{pmatrix}1 & 0\\0 & 1\end{pmatrix}, \ph C(t) =~ - K(t), \linebreak t\ge t_0$. Hence by Theorem 3.1 Eq. (3.28) is oscillatory provided is oscillatory the following scalar system
$$
\sist{\phi_1' = \psi_1;}{\psi_1' = - \mu(t) \phi_1, \ph t\ge t_0.}
$$
This system is equivalent to the second order scalar equation
$$
\phi_1'' + \mu(t) \phi_1 = 0, \phh t\ge t_0,
$$
which is oscillatory (see [15]). Therefore Eq. (3.28) is oscillatory.
It is not difficult to verify that the results of works [16 -20] are not applicable to Eq. (3.28).

Example 3.2.
Let
$$
B(t) = \begin{pmatrix} 1 & 1\\ 1 & 1\end{pmatrix}, \ph t\ge t_0. \eqno (3,29)
$$
 Then $\sqrt{B(t)} = \frac{\sqrt{2}}{2} \begin{pmatrix} 1 & 1\\ 1 & 1\end{pmatrix}, \ph \sqrt{B(t)}' \equiv~ 0,  \ph t\ge t_0,$ and  $ F(t) = \sqrt{2}\begin{pmatrix} 1 & 0\\ 0 & 1\end{pmatrix}, \ph t\ge t_0,$ is a solution of Eq. (3.21), on $[t_0;+\infty)$,
$$
P(t) = \begin{pmatrix} a_{11}(t) + a_{12}(t) & a_{11}(t) + a_{12}(t) \\ a_{21}(t) + a_{22}(t) &  a_{21}(t) + a_{22}(t) \end{pmatrix}, \phantom{aaaaaaaaaaaaaaaaaaaaaaaaaaaaaa} \eqno (3.30)
$$
$$
\phantom{aaaaaaaaaaaaaaaaaaaaa} Q(t) = (c_{11}(t) + 2 Re\hskip 2pt c_{12}(t) + c_{22}(t))B(t), \phh t\ge t_0.
 \eqno (3.31)
$$
Assume
$$
a_{11}(t) + a_{12}(t) = a_{21}(t) + a_{22}(t) \equiv 0, \phh t\ge t_0. \eqno (3.32)
$$
Then taking into account (3.30) and (3.31) we have: $\widetilde{\chi}_1(t) = \widetilde{\chi}_2(t)= - c_{11}(t) - 2 Re\hskip2pt  c_{12}(t) - \linebreak  - c_{22}(t), \ph t\ge t_0.$ Therefore by Corollary 3.1 under the restrictions (3.29) and (3.32) the system (1.1) is oscillatory provided the scalar equation
$$
\phi_1''(t)- [c_{11}(t) + 2 Re \hskip2pt c_{12}(t) + c_{22}(t)] \phi_1(t) = 0, \phh t\ge t_0,
$$
is oscillatory.

Assume now:
$$
a_{11}(t) + a_{12}(t) = a_{21}(t) + a_{22}(t) = \frac{\alpha}{t}, \ph c_{11}(t) + 2 Re \hskip2pt c_{12}(t) + c_{22}(t) = \frac{\alpha - \alpha^2}{t^2},  \eqno (3.33)
$$
$0\le \alpha \le 1, \ph t\ge 1$. Then taking into account (3.30) and (3.31) it is not difficult to verify that $\widetilde{\chi}_3(t) = \widetilde{\chi}_4(t) = \frac{\alpha^2 - \alpha}{t^2} \le 0, \ph t\ge 1.$ Hence by Theorem 3.4 under the restrictions (3.29) and (3.33)  the system (1.1) is non oscillatory.

Let now we assume:

\noindent
$\alpha_1) \ph a_{11}(t) + a_{12}(t) = a_{21}(t) + a_{22}(t) > 0, \ph t\ge t_0;$

\noindent
$\alpha_2) \ph a_{11}(t) + a_{12}(t)$ is increasing and continuously differentiable on $[t_0;+ \infty)$;

\noindent
$\alpha_3) \ph \frac{|(a_{11}(t) + a_{12}(t))' + c_{11}(t) + 2 Re \hskip2pt c_{12}(t) + c_{22}(t)|}{a_{11}(t) + a_{12}(t)} \le \lambda = const, \ph t\ge t_0.$

\noindent
Then taking into account (3.30) and (3.31) it is not difficult to verify that $\widetilde{\chi}_3(t) \le \lambda - [c_{11}(t) + 2 Re \hskip2pt c_{12}(t) + c_{22}(t)], \ph  \widetilde{\chi}_4(t) \le \lambda - [c_{11}(t) + 2 Re \hskip 2pt c_{12}(t) + c_{22}(t)], \ph t\ge t_0.$ Therefore by virtue of Theorem 3.4 under the restrictions (3.29) and $\alpha_1) - \alpha_3)$ the system (1.1) is non oscillatory.

{\bf Remark 3.4.} {\it Since under the restriction (3.29) $det B(t) \equiv 0, \ph t\ge t_0,$ the results of works [1 -9] are not applicable to the system (1.1) with (3.29).}

\vskip 20 pt

\centerline{ \bf References}

\vskip 20pt

\noindent
1. L. Li, F. Meng and Z. Zheng, Oscillation Results Related to Integral Averaging Technique\linebreak \phantom{a} for Linear Hamiltonian Systems, Dynamic Systems and Applications 18 (2009), \ph \linebreak \phantom{a} pp. 725 - 736.

\noindent
2. F. Meng and  A. B. Mingarelli, Oscillation of Linear Hamiltonian Systems, Proc. Amer.\linebreak \phantom{a} Math. Soc. Vol. 131, Num. 3, 2002, pp. 897 - 904.

\noindent
3. Q. Yang, R. Mathsen and S. Zhu, Oscillation Theorems for Self-Adjoint Matrix \linebreak \phantom{a}   Hamiltonian
 Systems. J. Diff. Equ., 19 (2003), pp. 306 - 329.

\noindent
4. Z. Zheng and S. Zhu, Hartman Type Oscillatory Criteria for Linear Matrix Hamiltonian  \linebreak \phantom{a} Systems. Dynamic  Systems and Applications, 17 (2008), pp. 85 - 96.

\noindent
5. Z. Zheng, Linear transformation and oscillation criteria for Hamiltonian systems. \linebreak \phantom{a} J. Math. Anal. Appl., 332 (2007) 236 - 245.

\noindent
6. I. S. Kumary and S. Umamaheswaram, Oscillation Criteria for Linear Matrix \linebreak \phantom{a} Hamiltonian Systems, Journal of Differential Equations, 165, 174 - 198 (2000).

\noindent
7. Sh. Chen, Z. Zheng, Oscillation Criteria of Yan Type for Linear Hamiltonian Systems, \linebreak \phantom{a} Computers and Mathematics with Applications, 46 (2003), 855 - 862.

\noindent
8. Y. G. Sun, New oscillation criteria for linear matrix Hamiltonian systems. J. Math. \linebreak \phantom{a} Anal. Appl., 279 (2003) 651 - 658.

\noindent
9. K. I. Al - Dosary, H. Kh. Abdullah and D. Husein. Short note on oscillation of matrix \linebreak \phantom{a} hamiltonian systems. Yokohama Mathematical Journal, vol. 50, 2003.

\noindent
10. G. A. Grigorian, Oscillatory and Non Oscillatory Criteria for the Systems of Two \linebreak \phantom{aa}   Linear First Order Two by Two Dimensional Matrix Ordinary Differential Equations. \linebreak \phantom{aa}   Archivum  Mathematicum, Tomus 54 (2018), PP. 189 - 203.

\noindent
11. G. A. Grigorian.  On Two Comparison Tests for Second-Order Linear  Ordinary\linebreak \phantom{aa} Differential Equations (Russian) Differ. Uravn. 47 (2011), no. 9, 1225 - 1240; trans-\linebreak \phantom{aa} lation in Differ. Equ. 47 (2011), no. 9 1237 - 1252, 34C10.

\noindent
12. G. A. Grigorian, "Two Comparison Criteria for Scalar Riccati Equations with\linebreak \phantom{aa} Applications". Russian Mathematics (Iz. VUZ), 56, No. 11, 17 - 30 (2012).

\noindent
13. G. A. Grigorian, On the Stability of Systems of Two First - Order Linear Ordinary\linebreak \phantom{aa} Differential Equations, Differ. Uravn., 2015, vol. 51, no. 3, pp. 283 - 292.

\noindent
14.  G. A. Grigorian. Oscillatory Criteria for the Systems of Two First - Order Linear\linebreak \phantom{a} Ordinary Differential Equations. Rocky Mountain Journal of Mathematics, vol. 47,\linebreak \phantom{a} Num. 5, 2017, pp. 1497 - 1524

\noindent
15. G. A. Grigorian, On one Oscillatory Criterion for The Second Order Linear
 Ordinary \linebreak \phantom{a} Differential Equations. Opuscula Math. 36, Num. 5 (2016), 589–601. \\ \phantom{a}
   http://dx.doi.org/10.7494/OpMath.2016.36.5.589

\noindent
16. L. H. Erbe, Q. Kong and Sh. Ruan, Kamenev Type Theorems for Second Order Matrix\linebreak \phantom{aa}  Differential Systems. Proc. Amer. Math. Soc. Vol. 117, Num. 4, 1993, 957 - 962.

\noindent
17. R. Byers, B. J. Harris and M. K. Kwong, Weighted Means and Oscillation Conditions\linebreak \phantom{a}  for Second Order Matrix Differential Equations. Journal of Differential Equations\linebreak \phantom{a} 61, 164 - 177 (1986).

\noindent
18. G. J. Butler, L. H. Erbe and A. B. Mingarelli, Riccati Techniques and Variational\linebreak \phantom{aa} Principles in Oscillation Theory for Linear Systems, Trans. Amer. Math. Soc. Vol. 303,\linebreak \phantom{aa} Num. 1, 1987, 263 - 282.

\noindent
19. A. B. Mingarelli, On a Conjecture for Oscillation of Second Order Ordinary Differential\linebreak \phantom{aa} Systems, Proc. Amer. Math. Soc., Vol. 82. Num. 4, 1981, 593 - 598.

\noindent
20. Q. Wang, Oscillation Criteria for Second Order Matrix Differential Systems Proc.\linebreak \phantom{aa} Amer.  Math. Soc. Vol. 131, Num. 3, 2002, 897 - 904.

\end{document}